\documentclass[reqno]{amsart}
\usepackage{url}
\usepackage{amssymb}
\usepackage{graphicx}
\usepackage[colorinlistoftodos]{todonotes}
\usepackage{amsmath}
\usepackage{amsthm}

\newtheorem{theorem}{Theorem}[section]
\newtheorem{lemma}[theorem]{Lemma}

\newtheorem{conjecture}[theorem]{Conjecture}

\theoremstyle{definition}

\theoremstyle{remark}

\numberwithin{equation}{section}
\usepackage{algorithmic}
\usepackage[top=3cm,bottom=2cm,right=2cm,left=2cm]{geometry}

\begin{document}


\title[On a problem of Gordon and Hughes]{On a problem of Gordon and Hughes:\\ The 11-adic properties of meromorphic modular functions on $\mathrm{X}_0(22)$}


\author{Koustav Banerjee}
\address{}
\curraddr{}
\email{}
\thanks{}

\author{Nicolas Allen Smoot}
\address{}
\curraddr{}
\email{}
\thanks{}



\date{}

\dedicatory{Dedicated to Krishnaswami Alladi, in honor of his 70th birthday}

\begin{abstract}
Congruence families occur throughout the theory of modular forms.  These families superficially resemble one another, though they can vary substantially in their overall difficulty.  One major complication resides in the topology of the modular curve associated with a given congruence family.  All previously studied families are associated with curves of genus 0 or 1.  The existence of a family associated with a curve of genus $2$ (and the associated difficulties of a proof thereof) was first speculated upon by Gordon and Hughes.  We give the first explicit congruence family associated with $\mathrm{X}_0(22)$, a modular curve of genus $2$, and prove it.  We also propose two families of conjectures associated with the same modular curve.
\end{abstract}

\maketitle

\section{Introduction}

In 1981 Gordon and Hughes studied Ramanujan type congruence families for the distinct partition function, denoted by $q(n)$, modulo powers of $5$ and $7$. They \cite[Section 1]{Gordon} speculated on the possibility of analogous congruences by powers of $11$: 
\begin{quote}
\textit{``As far as congruences for $q(n)$ mod $11^{\gamma}$ are concerned, the problem is greatly complicated by the fact that the group $\Gamma_0(22)$ has genus 2.  The computational difficulties seem too formidable to justify pursuing this problem here."}
\end{quote}

It is nearly a half-century since this paper was published, and no congruence families have been found associated with $\Gamma_0(22)$ or any other congruence subgroup corresponding to a genus $2$ modular curve. In this paper, we give an explicit example of such a family and propose a broader class of such genus $2$ congruence families. The congruences are exhibited by the $2$-elongated plane partition diamond function $d_2(n)$ rather than $q(n)$; nevertheless, the associated group is $\Gamma_0(22)$, as Gordon and Hughes had anticipated. Likewise, the difficulties, theoretical and computational, are as formidable as Gordon and Hughes foresaw. The principal theoretical novelty of the proof is the use of localizing functions that are not modular units, i.e., whose divisors are not supported on the cusps. We believe that the methods developed to handle such functions have implications beyond the present family, both for modular congruences and for the $\ell$-adic theory of modular forms more generally.

\subsection{Background}

The theory of modular congruence relationships has been an active area of research since Ramanujan's groundbreaking work \cite{Ramanujan} on the divisibility properties of the partition function. Consider the original congruence families discovered and proved by Ramanujan, Watson \cite{Watson}, and Atkin \cite{Atkin}:
\begin{align}
\text{If }\ n,\alpha\in\mathbb{N}\  \text{ and }\  &24n\equiv 1\ \left(\textnormal{mod}\ \ell^{\alpha}\right),\text{ then } p(n)\equiv 0\ \left(\textnormal{mod}\ \ell^{\beta}\right),\label{origconjecR5711}\\
\text{with }\ \ &\beta := 
\begin{cases}
\alpha & \text{ if }\ \ell\in\{5,11\},\\
\left\lfloor \frac{\alpha}{2} \right\rfloor + 1 & \text{ if }\ \ell=7.\nonumber
\end{cases}
\end{align}
Since the proof of these families, the arithmetic of the partition function has been an active area of research. In his seminal work, Ono~\cite{Ono2000} proved that for every prime $\ell\ge5$ there are infinitely many congruences of the form $p(An+B)\equiv0\pmod{\ell}$, and Ahlgren~\cite{Ahlgren2000} extended this to all moduli coprime to $6$ (see also Ahlgren and Ono~\cite{AhlgrenOno2001}). Ahlgren and Boylan \cite{AhlgrenBoylan2003} used tools from the theory of Galois representations and modular forms to show that Ramanujan's partition congruences are the only ones which are of the form $p(\ell n+\beta)\equiv0\pmod{\ell}$. For prime powers, Folsom, Kent, and Ono \cite{FKO} and Boylan and Webb \cite{BW} established a strong form of $\ell$-adic stability for the $U_\ell$-action. Ahlgren, Allen, and Tang \cite{AhlgrenAllenTang2022} proved that congruences like Atkin's hold modulo every prime $\ell\ge 5$.

This richness of arithmetic properties has been found to exist within the coefficients of many different modular forms. We are particularly interested in $\ell$-adic patterns analogous to (\ref{origconjecR5711}). A striking aspect of the subject is the substantial variability in the difficulty of a given congruence family. Some admit short proofs by methods that were already available to Ramanujan. Others remain open to this day. 

One major difficulty lies in the topology of the associated modular curve, e.g., the cusp count and the genus. The most difficult problems generally occur when the associated modular curve has many cusps. Naively speaking, in these cases the action of the relevant $U_\ell$ operator can no longer be controlled within a single, finitely generated space of modular functions, and the $\ell$-adic estimates required to close the induction break down. Families of this kind, such as those arising from traces of singular moduli and Hurwitz class numbers, have so far resisted even the most sophisticated techniques currently available \cite{Beazer, Guerzhoy}.

Another significant factor in the difficulty of a congruence family lies with the genus. For example, the associated modular curve for \eqref{origconjecR5711} is simply $\mathrm{X}_0(\ell)$.  In all three cases the cusp count is $2$.  However, for $\ell\in\{5,7\}$, the curve has genus $0$, while the curve has genus $1$ for $\ell=11$.  Watson, who had not completed the proof of \eqref{origconjecR5711} for $\ell=11$, speculated in 1938 that the proof would not be especially difficult, but  perhaps ``recht langweilig" \cite{Watson}.  Upon his completion of the proof for $\ell=11$ in 1967, Atkin \cite[Section 5]{Atkin} agreed with Watson's speculation. 

It has long been understood that families associated with curves of positive genus are difficult to treat by conventional techniques. In addition to the aforementioned remarks by Gordon and Hughes, Calegari \cite[Section 4.2]{Calegari} has also pointed out that the standard proof methods ``become essentially infeasible as soon as [the associated modular curve] has genus $>0$". The difficulty of a proof notwithstanding, it is interesting to note that no explicit congruence families for genus $2$ appear to have been conjectured prior to this work. To the best of our knowledge, the possibility of such a family has been raised in the literature only by Gordon and Hughes \cite{Gordon}. To test the expectation that such a family should be accessible to current methods, we undertook a search for a suitable example. A computational search for $q(n)$, the function originally considered by Gordon and Hughes, produced no such candidate. A natural alternative is provided by the Fourier coefficients of the following eta quotient, whose congruence families have already served as a valuable testing ground for genre $0$ and $1$ \cite{Banerjee1, Banerjee2}, and which therefore suggested itself as a promising source of a genus $2$ example:
\begin{align*}
D_k(q) := \sum_{n\ge 0}d_k(n)q^n = \frac{(q^2;q^2)_{\infty}^k}{(q;q)_{\infty}^{3k+1}}.
\end{align*}  The Fourier coefficients $d_k(n)$ have an elegant combinatorial interpretation, discovered by Andrews and Paule \cite{AndrewsPaule} in their ``Partition Analysis" series of papers. It is among the Fourier coefficients of $D_2(q)$ that we find a relevant family of interest, the main result of this paper:
\begin{theorem}\label{Thm12}
Let $n,\alpha\in\mathbb{Z}_{\ge 1}$ such that $8n\equiv 1\pmod{11^{\alpha}}$.  Then 
$$
d_2(n)\equiv 0\ \left({\textnormal{mod}\ 11^{\left\lfloor\frac{\alpha+1}2\right\rfloor}}\right).$$
\end{theorem}

\subsection{Localization and difficulty of the proof}

Gordon and Hughes’s assessment of the difficulty proves accurate. The genus of $\mathrm{X}_0(22)$, together with difficulties in its corresponding space of eta quotients, necessarily complicates our attempts at a proof of Theorem \ref{Thm12}. As in previous approaches (see \cite{Banerjee1,Banerjee2}) via modular equation manipulations, the idea was straightforward: to study $d_2(n)$ in the progressions relevant to us, we construct a sequence of functions $L_{\alpha}$ of the form
\begin{align*}
L_{\alpha} &:= \phi_{\alpha}\sum_{8n\equiv 1\bmod{11^{\alpha}}} d_2(n)q^{\left\lfloor \frac{n}{11^{\alpha}}\right\rfloor+1}.
\end{align*}  With the right holomorphic prefactors $\phi_{\alpha}$, we can arrange so that $L_{\alpha}$ is a modular function over $\Gamma_0(22)$, i.e., a meromorphic function on $\mathrm{X}_0(22)$ with poles only at the cusps.  We can construct a pair of modified Atkin $U_{11}$ operators $U^{(i)}$, with $i\in\{0,1\}$, such that
\begin{align*}
U^{(a)}\left( L_{\alpha} \right) = L_{\alpha+1},
\end{align*} with $a\equiv \alpha\bmod{2}$.
We want to place $L_{\alpha}$ into a suitable subspace of the meromorphic functions over $\mathrm{X}_0(22)$, in which the basis functions are subject to certain modular equations.  We then want to study how $U^{(i)}$ affects those basis functions.  This is done by directly examining a handful of initial cases, and then using the associated modular equations to construct more general relationships. Recall that $\mathrm{X}_0(22)$ is a compact Riemann surface with cusps
\begin{align*}
\left\{[\infty], \left[\frac{1}{11}\right], \left[\frac 12\right], [0]\right\}.
\end{align*} 
Since $\mathrm{X}_0(22)$ is compact, every non-constant meromorphic function on it has a pole.  Generally, we prefer to restrict these poles to the cusps. We consider the space of holomorphic functions over the curve except for a pole (of possibly large order) at the cusp $[0]$, and denote this space as $\mathcal{M}^{(0)}( \mathrm{X}_0(22))$. Using a straightforward computational procedure, we can use the Dedekind eta quotient to construct suitable functions with integer Fourier coefficients, and with orders $-3$, $-4$, $-5$ at $[0]$.  We denote these functions $x, y, z$, respectively.  Their explicit form is given below (see \eqref{xdef}--\eqref{zdef}).

Combinations of these three functions can be used to construct functions of order $-n$ for $n\ge 6$.  The only missing orders are $-1$, $-2$.  By the Weierstrass gap theorem \cite{Paule2} exactly two pole orders at $[0]$
are not attained, i.e., as $x$, $y$, $z$ have orders $-3$, $-4$, $-5$, the gaps are $1$ and $2$. Thus, we have
\begin{align*}
\mathcal{M}^{(0)}\left( \mathrm{X}_0(22) \right)=\mathbb{C}[x]\oplus y\mathbb{C}[x]\oplus z\mathbb{C}[x].
\end{align*}  Suppose $L_{\alpha}$ has a pole at $[0]$ and nowhere else. Then we could express it as a suitable polynomial combination of $x,y$, and $z$.  However, the poles of $L_{\alpha}$ are supported at $[0]$ and $[\frac 12]$.  In context of the localization method, we can multiply $L_{\alpha}$ by some modular eta quotient which lives in $\mathcal{M}^{(0)}( \mathrm{X}_0(22))$, but has positive order at $[\frac 12]$.  Indeed, a very natural candidate exists. Consider,\footnote{The second identity can be verified easily using \eqref{ydef} and \eqref{zdef}.}
\begin{align*}
w := \frac{(q^2;q^2)_{\infty}^{11}(q^{11};q^{11})_{\infty}}{(q;q)_{\infty}^{11}(q^{22};q^{22})_{\infty}} = 1+11y+11^2z.
\end{align*}
Indeed, $w^3 L_1\in \mathcal{M}^{(0)}\left( \mathrm{X}_0(22) \right)$.  We can therefore express $L_1$ as a rational polynomial with summands of the form
\begin{align*}
\frac{y^jz^kx^m}{(1+11y+11^2z)^3}\ \text{ with}\ m\in \mathbb{N}_0,\ j,k\in\{0,1\},\ \text{and}\  j+k\le 1.
\end{align*} 
However, such an approach raises various difficulties upon attempts to apply it to $L_{\alpha}$ generally.  In the first place, the denominator is governed by a degree $55$ modular equation, while $x$ is governed by a degree $33$ modular equation.  As such, we will be forced to examine how $U^{(i)}$ affects $1815$ different combinations, with an additional factor of $2$ for the value of $i$, and a factor of $3$ for the possible presence of $y$ and $z$.  This immediately forces well over $10^4$ large polynomial relations to be verified. Moreover, the presence of $y$ and $z$ in $w$ ensures that one must tediously remove monomials $y^az^b$ with $a+b\ge 2$ from resultant representations. Replacing $z$ by $w$ does not resolve this complication, because such a substitution would incur high powers of $11$ in the denominators of coefficients in $L_{\alpha}$. This immediately destroys the utility of $L_{\alpha}$ as a witness to divisibility by powers of $11$. Another alternative is to work with a $\mathbb{C}[w,w^{-1}]$ module of rank $5$.  In this case, one still requires the localization technique and needs to verify at least $3\cdot 10^4$ associated relations. Finally, the methods of Paule and Radu \cite{Paule}, and of Garvan, Sellers, and the second author \cite[Section 3]{Garvan0} do not appear to adapt to the present setting.

These obstacles are consistent with Gordon and Hughes’s original assessment. We instead construct the localizing factor from $x$ alone. To remove the poles of $L_{\alpha}$ at $[\frac 12]$, note that the expansion of $x$ near $[\frac 12]$ in terms of the local variable $\hat{q}$ is the following
\begin{align*}
x\left( \frac{\tau}{2\tau+1} \right) = -\frac{2}{11} + b\hat{q}+O\left(\hat{q}^2\right)\ \ \text{with}\ \ \ b\neq 0.
\end{align*}  Thus ${\rm ord}_{[\frac 12]}(2+11x)=1$ and 
\begin{align}
(2+11x)^{\psi}L_{\alpha}\in\bigoplus_{v\in\{x,y,z\}} v\,\mathbb{C}[x]\label{Lalphainto0}
\end{align} for an appropriate positive integer power $\psi$ depending on $\alpha$. We can now apply the modular cusp analysis to directly check just $198$ expressions for $U^{(i)}\left(y^jz^kx^m\right)$.  With these relations, we move on to examine some additional $6534$ cases by algebraic manipulation.  This is substantially smaller than the other possibilities considered. As an example, we consider the very first case:
\begin{align*}
L_{1} &= \frac{(q^{11};q^{11})^{7}_{\infty}}{(q^{22};q^{22})^2_{\infty}} \sum_{n=0}^{\infty}d_2\left(11n + 7\right)q^{n+2}.
\end{align*}  Next we rewrite $L_1$ as\,\footnote{The explicit forms of $M_v$ are given in Section 1 of the Mathematica supplement. Details on the availability of the notebook are given in Subsection 1.3.}
\begin{align}
L_1 = \frac{11}{2(2+11x)^{15}}\left( M_x(x) + y M_y(x) + z M_z (x) + (3+7x) r_L \right)\label{L1identity11},
\end{align} with
\begin{align*}
M_v\in\mathbb{Z}[X],\ v\in\{x,y,z\},\quad r_L := \frac{1}{11} \left(x+3x^2+9x^3+10y+10xy+2x^2y+7z+8xz\right)\in\mathbb{Z}[[q]].
\end{align*}   The function $r_L$, which we call the anomalous term, plays an important role. Its numerator is divisible by $11$, although this is not visible from the individual coefficients. We prove that $r_L$ is indeed an integer power series by appealing to Sturm's theorem.

The more general result is given below. From this, Theorem \ref{Thm12} follows:
\begin{theorem}\label{theorem0}
Let $L_{\alpha}$ be given as (\ref{Lall}), and
\begin{align}
\psi :=& \psi(\alpha) = \left\lfloor \frac{11^{\alpha+1}}{8} \right\rfloor +1 -\mathrm{gcd}(\alpha,2),\label{psidefn}\text{ and }\beta := \beta(\alpha) = \left\lfloor \frac{\alpha+1}{2} \right\rfloor.
\end{align}  There exist sequences $(j_{\alpha})_{\alpha\ge 1},(k_{\alpha})_{\alpha\ge 1}\in \mathbb{Z}\left[\frac 12\right]$ such that for all $\alpha\ge 1$, $j_{2\alpha}=k_{2\alpha}=0$, and
\begin{align*}
\frac{(2+11x)^{\psi}}{11^{\beta}} L_{\alpha} + (j_{\alpha}+k_{\alpha}x)\,r_L \in\bigoplus_{v\in\{x,y,z\}} v\,\mathbb{Z}\left[\frac{1}{2},x\right].
\end{align*}
\end{theorem}
It is worthwhile to note that in \cite{Banerjee1, Banerjee2}, the analogous localizing factors are eta quotients, and therefore \textit{modular units}. But, in this case, $2+11x$ is not a modular unit. Consequently, we must work in a larger localization ring, i.e., functions of the form $y^jz^kx^m c(x)^{-d}(2+11x)^{-n}$ with $n,d\ge 0$, and a certain polynomial $c(x)$ (see Section \ref{Appen}) induced by the modular equations for $x$. The heart of the proof is to show that, when $L_{\alpha}$ is viewed in this ring, the factor $c(x)^d$ cancels without destroying the $11$-adic structure. This is the novelty and significance of the framework we propose. 

\subsection{Outline}

In Section \ref{constructingLalpha} we give an exact definition to $L_{\alpha}$, together with some of its key meromorphic properties.  We also define the operators $U^{(i)}$ and their effect on $L_{\alpha}$.  In Section \ref{section3} we precisely define $x, y, z$, as well as the modular equation for $x$, and the resulting recurrence relations needed to treat the associated function spaces.  In Section \ref{section4} we study the action of $U^{(i)}$ on the function spaces and $r_L$.  In Section \ref{cancelcx} we account for the additional factors induced by the non-unitary behavior of the associated localizing term.  We give the complete proof of Theorem \ref{theorem0} in Section \ref{completeproofA}.  In Section \ref{proofoffundamentalrelations} we show how to prove the 198 initial relations required for Theorem \ref{theoremcrit} in Section \ref{section3}, together with the expression for $L_1$ above, and the modular equation for $x$.  We end in Section \ref{sectconcl} with two families of conjectures also over $\mathrm{X}_0(22)$, and briefly discuss the relevance of proof techniques different from the modular equation approach. Section \ref{Appen} provides explicit descriptions of the constants in Theorem \ref{theoremcrit}, as well as the function $c(x)$ mentioned in Section \ref{section3}.

All the supplementary computations are documented in three Mathematica notebooks.\footnote{The Mathematica notebooks can be found at \url{https://homepage.univie.ac.at/nicolas.allen.smoot/?page_id=257}}. Moreover, we have also verified these computations with PARI/GP.\footnote{The reader can access the files from \url{https://sites.google.com/view/koustav-banerjee/research?authuser=0}.}

\section{The Function Space of $\mathrm{X}_0(22)$}\label{modularfunctionmanips}

\subsection{Construction of $L_{\alpha}$}\label{constructingLalpha}

To construct $L_{\alpha}$,  we define, for $\tau\in \mathbb{H}$,
\begin{align*}
\mathcal{A}&:= \mathcal{A}(\tau):= q^{15}\frac{D_2(q)}{D_2(q^{121})} = q^{15}\frac{(q^2;q^2)^2_{\infty}(q^{121};q^{121})^{7}_{\infty}}{(q;q)^{7}_{\infty}(q^{242};q^{242})^2_{\infty}},\\
U^{(1)}(f) :=& U_{11}(f),\quad U^{(0)}(f) := U_{11}(\mathcal{A}\cdot f),\quad U^{(\alpha)}(f) := U^{(i)}(f),\quad  \alpha\equiv i\pmod{2},
\end{align*}
 and the function sequence $\mathcal{L}:=\left( L_{\alpha} \right)_{\alpha\ge 1}$ by
\begin{align}
L_{\alpha}(\tau) &:= D_2\left(q^{11^{(\alpha)}}\right)^{-1}\cdot\sum_{8n\equiv 1\bmod{11^{\alpha}}} d_2(n)q^{\left\lfloor \frac{n}{11^{\alpha}} \right\rfloor+1+(\alpha)},\label{Lall}
\end{align} with $(\alpha)\in\{0,1\}$ such that $(\alpha)\equiv\alpha\pmod{2}$.  By convention, we set $L_0:=1$. 
\begin{lemma}\label{use}
For all $\alpha\ge 0$, we have
\begin{align*}
L_{\alpha+1} &= U^{(\alpha)}\left( L_{\alpha} \right)\in\mathcal{M}\left( \mathrm{X}_0(22) \right).
\end{align*}
\end{lemma}
\begin{proof}
This follows from a routine computation with eta quotients, see \cite[Lemma 4.1]{Banerjee1}, \cite[Lemma 2.7]{Banerjee2}, and \cite[Lemma 2.13]{Smoot2} for analogous arguments. 
\end{proof}
Since $\mathcal{A}$ is an eta quotient, it has no zeros or poles in $\mathbb{H}$, hence the poles of $L_{\alpha}$ lie at the cusps of $\mathrm{X}_0(22)$. As noted before,\footnote{Here, we use a different representative $[1]=[0]$ and $[\frac1{22}]=[\infty]$.} the cusps of $\mathcal{M}\left( \mathrm{X}_0(22) \right)$ are $[1]$, $[\frac 12]$, $[\frac{1}{11}]$, $[\frac{1}{22}]$.  We study the behavior of each $L_{\alpha}$ at these cusps.

From \eqref{Lall}, it is not hard to see that $L_{\alpha}$ has positive order at $[\frac{1}{22}]$.  For the remaining three cusps we use bounds given by \cite[Definition 38, Theorem 47]{Radu}.  Let us denote the order of a function at the cusp $[\frac{a}{c}]$ as $\mathrm{ord}_{\frac ac}$.  Then we have
\begin{align*}
\mathrm{ord}_{1}\left( L_{2\alpha-1} \right)&\ge\frac{-11^{2\alpha}+1}{2},\quad \mathrm{ord}_{1}\left( L_{2\alpha} \right) \ge\frac{-11^{2\alpha+1}+11}{2},\quad \mathrm{ord}_{\frac{1}{11}}\left( L_{2\alpha-1} \right) \ge\frac{11^{2\alpha-2}-1}{2\cdot 11^{2\alpha-1}} \ge 0,\\
\mathrm{ord}_{\frac{1}{11}}\left( L_{2\alpha} \right)&\ge\frac{11^{2\alpha}-1}{2\cdot 11^{2\alpha}} \ge 0,\end{align*} and 
\begin{align}
\mathrm{ord}_{\frac 12}\left( L_{2\alpha-1} \right) \ge\frac{-11^{2\alpha}+1}{8} = -\psi(2\alpha-1),\qquad \mathrm{ord}_{\frac 12}\left( L_{2\alpha} \right) \ge\frac{-11^{2\alpha+1}+11}{8} = -\psi(2\alpha).\label{Lalphaordersat12}
\end{align}
This proves that $L_{\alpha}$ has negative order only at $[0]$ and $[\frac 12]$, with well-defined lower bounds in both cases.  Moreover, (\ref{Lalphaordersat12}), together with the positive order of $2+11x$ at $[\frac{1}{2}]$, confirms (\ref{Lalphainto0}).

\subsection{Basis Functions For $\mathrm{X}_0(22)$}\label{section3}

We come now to the question of constructing a suitable set of basis functions with which to express $L_{\alpha}$.  Related functions $x,y,z$ were already computed with a pole at $[\infty]$ in \cite[Section 4.1]{Radu}, but we quickly explain the concept.  We want functions of orders $-3$, $-4$ and $-5$ at $[0]$. The space of modular eta quotients over $\Gamma_0(22)$ with a pole only at $[0]$ is a monoid which is generated by $7$ elements.  These can be quickly computed with existing algorithmic machinery.  Of these elements, $3$ of them have order $-5$ at $[0]$, while one has order $-6$, one has order $-7$, one with order $-8$, and one with order $-9$. Thus we can easily find a suitable function $z$.  To get $x$ and $y$ we can simply compute the expansions of different combinations of order $-5$ eta quotients at $[0]$ and subtract off the higher order terms.  This can be done in different ways, but one suitable example is the following:
\begin{align}
x &= \frac{1}{44}\left( -7 + 11\frac{(q^{2};q^{2})^{4}_{\infty}(q^{11};q^{11})^{8}_{\infty}}{(q;q)^{8}_{\infty}(q^{22};q^{22})^{4}_{\infty}} - 4\frac{(q^{2};q^{2})^{11}_{\infty}(q^{11};q^{11})_{\infty}}{(q;q)^{11}_{\infty}(q^{22};q^{22})_{\infty}} + 132q^5 \frac{(q^{2};q^{2})^{3}_{\infty}(q^{22};q^{22})^{7}_{\infty}}{(q;q)^{7}_{\infty}(q^{11};q^{11})^{3}_{\infty}} \right),\label{xdef}\\
y &= \frac{1}{11}\left( -1 + \frac{(q^{2};q^{2})^{11}_{\infty}(q^{11};q^{11})_{\infty}}{(q;q)^{11}_{\infty}(q^{22};q^{22})_{\infty}} - 121q^5 \frac{(q^{2};q^{2})^{3}_{\infty}(q^{22};q^{22})^{7}_{\infty}}{(q;q)^{7}_{\infty}(q^{11};q^{11})^{3}_{\infty}} \right),\label{ydef}\\
z &= q^5 \frac{(q^{2};q^{2})^{3}_{\infty}(q^{22};q^{22})^{7}_{\infty}}{(q;q)^{7}_{\infty}(q^{11};q^{11})^{3}_{\infty}}.\label{zdef}
\end{align}  
Since $(q;q)^{11}_n\equiv (q^{11};q^{11})_n\pmod{11}$, the Fourier coefficients of $x$ and $y$ are integers. Note that $r_L$ is also a series with integer coefficients, which is verified by checking the first $361$ coefficients (Sturm bound) of the weight $120$ form $f_r:=\eta(22\tau)^{240}r_L$.\footnote{The verification is done with Mathematica, see in Section 2b of the Mathematica supplement.}

We now derive a recurrence for the action of $U^{(i)}$ on powers of $x$. Both $x=x(\tau)$ and $x(11\tau)$ are modular functions for $\Gamma_0(242)$. On $\mathrm{X}_0(242)$, the function $x$ has a single pole, of order $33$, at $[0]$, whereas $x(11\tau)$ has poles only at $[0]$ and $[\frac1{11}]$, of total order $33$. Hence, by \cite[Theorem 8.3]{Forster}, $x$ satisfies a monic polynomial equation of degree $33$ over $\mathbb{C}(x(11\tau))$ . Since every pole of $x$ is also a pole of $x(11\tau)$, the function $x$ is integral over $\mathbb{C}[x(11\tau)]$, so the coefficients of this equation are polynomials: there exist $a_j\in\mathbb{C}[X]$, $0\le j\le 32$, such that
\begin{equation}\label{modX}
x^{33}+\sum_{j=0}^{32}a_j\bigl(x(11\tau)\bigr)\,x^j=0.
\end{equation}
In fact $a_j\in\mathbb{Z}[X]$.\footnote{The $a_j$ are too large to reproduce here and they are given in Section 2a of the supplement.} This is verified directly by multiplying the left-hand side of (3.4) by a suitable eta quotient that moves all poles to $[\infty]$, and check that the resulting principal part and constant term vanish (see Section 6).

For a polynomial $P$ we write $P(\tau):=P(x(\tau))$ and thus $U^{(i)}\bigl(P(11\tau)f\bigr)=P(\tau)\,U^{(i)}(f)$, which is how $a_\mu(\tau)$ and $b_\nu(\tau)$ arise in Lemma \ref{Lem1} below. Setting $t:=2+11x$, by \eqref{modX} with $x=\frac{t-2}{11}$, we get an associated modular equation for $2+11x$ that reads
\begin{align}
t^{33}+\sum_{k=0}^{32} b_k\left(x\left(11\tau\right)\right) t^k = 0.\label{modT}
\end{align}  Here $b_k\in\mathbb{Z}\left[\frac 12,x\right]$.\footnote{Explicit descriptions of $b_k$ are given in Section 2c of the Mathematica supplement.}  Precisely,\,\footnote{The expression of $c(x)$ is given in Section \ref{Appen}.}
\begin{align*}
b_0(\tau) &= -\frac{1}{2^{20}}(2+11x)^{11}c(x)\ \ \text{ with}\ \ c(x):= C_1(0) + \sum_{r=1}^{22} C_1(r) 11^{\left\lfloor\frac{11r+27}{13}\right\rfloor} x^r\in\mathbb{Z}[x],\\
C_1(0) &= 69550975755977284911104\equiv 4\pmod{11}.
\end{align*}
 
Given the form of $L_{\alpha}$ in Theorem \ref{theorem0}, and the linearity of $U^{(i)}$, we want to study
\begin{align}
U^{(i)}\left(\frac{y^jz^kx^m}{(2+11x)^n}\right).\label{uappliedtomonomial}
\end{align} 
Precisely, we have the following:
\begin{lemma}\label{Lem1}
For all $m,n\in\mathbb{Z}$, and $i,j,k\in\{0,1\}$, with $j+k\le 1$, we have
\begin{align*}
U^{(i)}\left( \frac{y^j z^k x^m}{(2+11x)^n} \right) = \frac{-2^{20}}{(2+11x)^{11} c(x)}\sum_{\mu=0}^{32}\sum_{\nu=1}^{33} a_{\mu}(\tau)b_{\nu}(\tau)\, U^{(i)}\left( \frac{y^j z^k x^{m+\mu-33}}{(2+11x)^{n-\nu}} \right).
\end{align*}
\end{lemma}
\begin{proof}
This is immediate from \eqref{modX} and \eqref{modT}.
\end{proof}
Due to the additional nontrivial factor $c(x)$, we study more general rational polynomial spaces of the following form:
\begin{align*}
\mathcal{W}^{(i)}_{n,d}:=& \left\{ \frac{1}{(2+11x)^{n} c(x)^d}\sum_{j+k\le 1}y^j z^k\sum_{m\ge 1-(j+k)} s_{jk}(m)\,11^{\theta^{i}_{jk}(m)} x^m \right\}.
\end{align*} 
Set
\begin{align*}
\theta^{i}_{jk}(m) := \mathrm{max}\left\{0,\left\lfloor \frac{11(m-4i+j+2k)+3i-7}{13} \right\rfloor\right\}.
\end{align*}  We now give a more useful recurrence:

\begin{theorem}\label{theoremcrit}
For each $(i,j,k,\gamma,\delta,m,n)\in\mathbb{N}^{7}_{0}$ with $i,j,k,\gamma,\delta\in\{0,1\}$, such that $j+k\le 1$ and $\gamma+\delta\le 1$, there exists $h^{ijk}_{\gamma\delta}(m,n,r)\in \mathbb{Z}\left[\frac 12\right]$ with finite support in $r$, such that
\begin{align}
U^{(i)}&\left( \frac{y^j z^k x^m}{(2+11x)^n} \right)\nonumber\\ =& \frac{1}{(2+11x)^{11n+15(1-i)} c(x)^{n}}\Bigg(k(1-i)6^m 4^n r_L + \sum_{0\le \gamma+\delta\le 1}\sum_{r\ge 1-(\gamma+\delta)} h^{ijk}_{\gamma\delta}(m,n,r)11^{\pi^{ijk}_{\gamma\delta}(m,r)}y^{\gamma} z^{\delta} x^r\Bigg), \label{Uyzxtotal}
\end{align} with
\begin{align*}
\pi^{ijk}_{\gamma\delta}(m,r) :&= \mathrm{max}\left\{0,\left\lfloor \frac{11r - m + \lambda^{ijk}_{\gamma\delta}}{13} \right\rfloor\right\},
\end{align*} and $\lambda^{ijk}_{\gamma\delta}\in\mathbb{Z}$ given in Section \ref{Appen}.
\end{theorem}
\begin{proof}
We prove by induction.  Suppose that \eqref{Uyzxtotal} holds for $m-\mu-1$ and $n-\nu$ with $m, n\in \mathbb{N}_{\ge 33}$ fixed and $0\le\mu\le 32$, $1\le\nu\le 33$.  We show that \eqref{Uyzxtotal} holds for $m,n\in \mathbb{N}_0$. By Lemma \ref{Lem1}, we have
\begin{align*}
&U^{(i)}\left( \frac{y^j z^k x^m}{(2+11x)^n} \right)= \frac{-2^{20}}{(2+11x)^{11} c(x)}\sum_{\mu=0}^{32}\sum_{\nu=1}^{33} a_{\mu}(\tau)b_{\nu}(\tau)\cdot U^{(i)}\left( \frac{y^jz^kx^{m+\mu-33}}{(2+11x)^{n-\nu}} \right)\\
&\hspace{1 cm}=\frac{-2^{20}}{(2+11x)^{11} c(x)}\sum_{\mu=0}^{32}\sum_{\nu=1}^{33} a_{\mu}(\tau)b_{\nu}(\tau) \frac{1}{(2+11x)^{11(n-\nu)+15(1-i)} c(x)^{n-\nu}}\\
&\hspace{2 cm}\times\left( k(1-i)6^{m+\mu-33}4^{n-\nu}r_L + \sum_{0\le \gamma+\delta\le 1}\sum_{r\ge 0} h^{ijk}_{\gamma\delta}(m+\mu-33,n-\nu,r)11^{\pi^{ijk}_{\gamma\delta}(m+\mu-33,r)}y^{\gamma} z^{\delta} x^r\right)\\
&\hspace{1 cm}=\frac{1}{(2+11x)^{11n+15(1-i)} c(x)^{n}}\left(\sum_{\substack{0\le\mu\le 32,\\ 1\le\nu\le 33}} W_{\mu,\nu}(\tau) 6^{m+\mu-33}4^{n-\nu}k(1-i) r_L\right)\\
&\hspace{2 cm}+\frac{1}{(2+11x)^{11n+15(1-i)} c(x)^{n}}\sum_{\substack{0\le\mu\le 32,\\ 1\le\nu\le 33,\\ 0\le \gamma+\delta\le 1,\\ r\ge 0}}W_{\mu,\nu}(\tau)h^{ijk}_{\gamma\delta}(m+\mu-33,n-\nu,r)11^{\pi^{ijk}_{\gamma\delta}(m+\mu-33,r)}y^{\gamma} z^{\delta} x^r,
\end{align*}  
where\,\footnote{Explicit descriptions of $w_{\mu,\nu}(\ell)$ is an integer-valued function with finite support in $\ell$ are given in Section 2e of the Mathematica supplement.}
\begin{align*}
W_{\mu,\nu}(\tau) :&= -2^{20}a_{\mu}(\tau)b_{\nu}(\tau)(2+11x)^{11(\nu-1)}c(x)^{\nu-1} =: \sum_{\ell\ge 0}w_{\mu,\nu}(\ell)11^{\phi(\ell,\mu)} x^{\ell},\\
\phi(\ell,\mu) :&= \mathrm{max}\left\{0,\left\lfloor \frac{11\ell+\mu-21}{13} \right\rfloor\right\}.
\end{align*}
Thus, the above becomes
\begin{align}
=& \frac{1}{(2+11x)^{11n+15(1-i)} c(x)^{n}}\sum_{\substack{0\le\mu\le 32,\\ 1\le\nu\le 33}} W_{\mu,\nu}(\tau) 6^{m+\mu-33}4^{n-\nu}k(1-i) r_L\nonumber\\
&+\frac{1}{(2+11x)^{11n+15(1-i)} c(x)^{n}}\sum_{\substack{0\le\mu\le 32,\\ 1\le\nu\le 33,\\ 0\le \gamma+\delta\le 1,\\ r\ge 0,\\ \ell\ge 0}}w_{\mu,\nu}(\ell)h^{ijk}_{\gamma\delta}(m+\mu-33,n-\nu,r)11^{\pi^{ijk}_{\gamma\delta}(m+\mu-33,r)+\phi(\ell,\mu)}y^{\gamma} z^{\delta} x^{r+\ell}.\label{u0mainvalterm}
\end{align} 

Next, we bound
\begin{align*}
\pi^{ijk}_{\gamma\delta}(m+\mu-33,r)+\phi(\ell,\mu) &\ge \left\lfloor \frac{11r - (m+\mu-33) + \lambda^{ijk}_{\gamma\delta}}{13} \right\rfloor + \left\lfloor \frac{11\ell+\mu-21}{13} \right\rfloor\\
&\ge \left\lfloor \frac{11(r+\ell) - m + \lambda^{ijk}_{\gamma\delta}+33  -21 -12}{13} \right\rfloor
\ge \pi^{ijk}_{\gamma\delta}(m,r+\ell).
\end{align*}  This confirms that the second sum in \eqref{u0mainvalterm}) has proper $11$-adic valuation.  Moreover, since $\ell\ge 0$, we have
\begin{align*}
r+\ell\ge r\ge 1-(\gamma+\delta).
\end{align*}  Indeed, the induction is complete except for the case that $i=0$ and $k=1$, in which case the first sum of \eqref{u0mainvalterm} remains to be examined. For the rest of the proof, assume $i=0$, $k=1$. Precisely, we now study the following function:
\begin{align}
W_{\mu,\nu}(\tau) 6^{m+\mu-33}4^{n-\nu} r_L &= -2^{20}6^{m+\mu-33}4^{n-\nu} a_{\mu}(\tau)b_{\nu}(\tau)(2+11x)^{11(\nu-1)}c(x)^{\nu-1}r_L.\label{wmunu6petc}
\end{align}

We show\,\footnote{Detail computation are given in the supplement notebook (see, Section 2f).} that when at least one of $\mu$, $\nu$ is not divisible by 11, (\ref{wmunu6petc}) expands into a polynomial with an 11-adic lower bound of
\begin{align*}
\mathrm{max}\left\{0,\left\lfloor \frac{11r-33+\lambda^{001}_{\gamma\delta}}{13} \right\rfloor\right\}\ge \mathrm{max}\left\{0,\left\lfloor \frac{11r-m+\lambda^{001}_{\gamma\delta}}{13} \right\rfloor\right\}=\pi^{001}_{\gamma\delta}(m,r)
\end{align*} 
for $y^{\gamma}z^{\delta}x^r$ (since $m\ge 33$).

The exceptional cases occur when both $\mu$ and $\nu$ are multiples of 11, i.e., when the pair $(\mu,\nu)=(11\lambda,11\kappa)$, for $0\le\lambda\le 2$, $1\le\kappa\le 3$. Combining these terms, we find\,\footnote{See Section 2f of the supplement file for detail calculations.} that 
\begin{align*}
\sum_{\substack{0\le\mu\le 32,\\ 1\le\nu\le 33}}W_{\mu,\nu}(\tau) 6^{m+\mu-33}4^{n-\nu} r_L &= f_0 + 6^{m-33}4^{n}\sum_{\substack{0\le\mu\le 2,\\ 1\le\nu\le 3}} W_{11\mu,11\nu}(\tau) 6^{11\mu}4^{-11\nu} r_L\\
&= f_0 + 6^{m-33}4^{n}\left( f_1 + \frac{1}{2^{24}}\left( f_2 + 2 r_L \right)\right),
\end{align*} where $f_0$, $f_1$, and $f_2$ satisfy the valuation requirements, and can be absorbed into \eqref{u0mainvalterm}. Thus, we next study $\frac{6^{m-33}4^{n}}{2^{23}} r_L$. Using the fact that
$2^{-23}6^{m-33}4^{n}\equiv 6^m 4^n\pmod{11}$, we have
\begin{align*}
\frac{6^{m-33}4^{n}}{2^{23}} r_L = 6^m 4^n r_L + f_3\ \ \text{with}\ \ f_3\in 11r_L\mathbb{Z}\left[\frac 12\right],
\end{align*} where $f_3$ again has the required $11$-adic valuation properties. This concludes\,\footnote{We give the computational details for the proof of the initial $6534$ cases is given in Section 2d of the supplementary file.} the proof of \eqref{Uyzxtotal}.
\end{proof}

Next we study the denominator of \eqref{uappliedtomonomial} by only varying its powers.  
\begin{lemma}\label{auxiliarycongruence}
Let $i,j,k\in\{0,1\}$ such that $j+k\le 1$ and $(i,k)\neq (0,1)$.  Let $m,r\in \mathbb{N}_0$, and $n\in \mathbb{N}_{\ge 33}$. We have
\begin{align*}
h^{ijk}_{\gamma\delta}(m,n,r)\equiv h^{ijk}_{\gamma\delta}(m,n-11,r) + 7h^{ijk}_{\gamma\delta}(m,n-22,r) + 9h^{ijk}_{\gamma\delta}(m,n-33,r)\pmod{11}.
\end{align*}
\end{lemma}
\begin{proof}
For $m$ fixed, by \eqref{modT}, we have
\begin{align}\label{auxiliarylinefirst}
&U^{(i)}\left(\frac{y^jz^kx^m}{(2+11x)^n}\right) = -\frac{1}{b_{0}}\sum_{\nu=1}^{33} b_{\nu} U^{(i)}\left(\frac{y^jz^kx^m}{(2+11x)^{n-{\nu}}}\right) =\frac{2^{20}}{(2+11x)^{11}c(x)}\sum_{\nu=1}^{33} b_{\nu} U^{(i)}\left(\frac{y^jz^kx^m}{(2+11x)^{n-{\nu}}}\right)\nonumber\\
&=\frac{2^{20}}{(2+11x)^{11}c(x)}\sum_{\nu=1}^{33} \frac{b_{\nu}}{(2+11x)^{11(n-\nu)+15(1-i)} c(x)^{n-\nu}}\sum_{\substack{0\le \gamma+\delta\le 1,\\r\ge 1-(\gamma+\delta)}} h^{ijk}_{\gamma\delta}(m,n-\nu,r)11^{\pi^{ijk}_{\gamma\delta}(m,r)}y^{\gamma} z^{\delta} x^r\nonumber\\
&= \frac{1}{(2+11x)^{11n+15(1-i)} c(x)^{n}}\sum_{\nu=1}^{33}2^{20}b_{\nu}(2+11x)^{11(\nu-1)}c(x)^{\nu-1}\sum_{\substack{0\le \gamma+\delta\le 1,\\r\ge 1-(\gamma+\delta)}} h^{ijk}_{\gamma\delta}(m,n-\nu,r)11^{\pi^{ijk}_{\gamma\delta}(m,r)}y^{\gamma} z^{\delta} x^r.
\end{align}  Define
\begin{align*}
Q(\nu,\tau) :&= 2^{20}b_{\nu}(2+11x)^{11(\nu-1)}c(x)^{\nu-1} =: \sum_{\ell\ge 0}Q_{\nu}(\ell)11^{\kappa(\ell,\nu)} x^\ell,
\end{align*}
with\,\footnote{Explicit descriptions are given in Section 2g of the supplement file.}
\begin{align*}
\kappa(\ell,\nu) :&= 
\begin{cases}
0 & \nu\in\{11,22,33\}\text{ and } \ell=0,\\
\left\lfloor \frac{11\ell+25}{13} \right\rfloor & \text{otherwise}.
\end{cases}
\end{align*}  
Thus, \eqref{auxiliarylinefirst} becomes
\begin{align}
&U^{(i)}\left(\frac{y^jz^kx^m}{(2+11x)^n}\right)\nonumber\\
&= \frac{1}{(2+11x)^{11n+15(1-i)} c(x)^{n}}\sum_{\substack{1\le\nu\le 33,\\ \ell\ge 0}}Q_{\nu}(\ell)11^{\kappa(\ell,\nu)} x^\ell\sum_{\substack{0\le \gamma+\delta\le 1,\\r\ge 1-(\gamma+\delta)}} h^{ijk}_{\gamma\delta}(m,n-\nu,r)11^{\pi^{ijk}_{\gamma\delta}(m,r)}y^{\gamma} z^{\delta} x^r\nonumber\\
&= \frac{1}{(2+11x)^{11n+15(1-i)} c(x)^{n}}\sum_{\substack{1\le\nu\le 33,\\ \ell\ge 0,\\ 0\le \gamma+\delta\le 1,\\r\ge 1-(\gamma+\delta)}} Q_{\nu}(\ell) h^{ijk}_{\gamma\delta}(m,n-\nu,r)11^{\pi^{ijk}_{\gamma\delta}(m,r)+\kappa(\ell,\nu)}y^{\gamma} z^{\delta} x^{r+\ell}.\label{auxcoeffmodt}
\end{align} 
On the other hand, since $(i,k)\neq(0,1)$, the anomalous term in \eqref{Uyzxtotal} vanishes, so \eqref{Uyzxtotal} expresses the left-hand side of \eqref{auxcoeffmodt} in terms of $h^{ijk}_{\gamma\delta}(m,n,r)$. Thus, we obtain
\begin{align}
&U^{(i)}\left(\frac{y^jz^kx^m}{(2+11x)^n}\right) = \frac{1}{(2+11x)^{11n+15(1-i)} c(x)^{n}}\sum_{\substack{0\le \gamma+\delta\le 1,\\r\ge 1-(\gamma+\delta)}} h^{ijk}_{\gamma\delta}(m,n,r)11^{\pi^{ijk}_{\gamma\delta}(m,r)}y^{\gamma} z^{\delta} x^r.\label{auxcoefforig}
\end{align}  As such, \eqref{auxcoeffmodt} and \eqref{auxcoefforig} give us a means to study the auxiliary coefficients $h^{ijk}_{\gamma\delta}(m,n,r)$.

Comparing the coefficients of $y^\gamma z^\delta x^r$ in \eqref{Uyzxtotal} and \eqref{auxcoeffmodt} gives
\[
h^{ijk}_{\gamma\delta}(m,n,r)\,11^{\pi^{ijk}_{\gamma\delta}(m,r)}=\sum_{\nu=1}^{33}\sum_{\ell=0}^{r}Q_\nu(\ell)\,h^{ijk}_{\gamma\delta}(m,n-\nu,r-\ell)\,11^{\pi^{ijk}_{\gamma\delta}(m,r-\ell)+\kappa(\ell,\nu)}.
\]
Unless $\ell=0$ and $\nu\in\{11,22,33\}$, we have \[\pi^{ijk}_{\gamma\delta}(m,r-\ell)+\kappa(\ell,\nu)\ge\pi^{ijk}_{\gamma\delta}(m,r)+1.\] Dividing by $11^{\pi^{ijk}_{\gamma\delta}(m,r)}$ and reducing modulo $11$, we obtain
\begin{align*}
h^{ijk}_{\gamma\delta}(m,n,r)\equiv{}&Q_{11}(0)h^{ijk}_{\gamma\delta}(m,n-11,r)+Q_{22}(0)h^{ijk}_{\gamma\delta}(m,n-22,r)\\
&+Q_{33}(0)h^{ijk}_{\gamma\delta}(m,n-33,r)\pmod{11},
\end{align*}
and the lemma follows from $Q_{11}(0)\equiv1$, $Q_{22}(0)\equiv7$, $Q_{33}(0)\equiv9\pmod{11}$.
\end{proof}

\section{Action of the $U^{(i)}$ Operators}\label{section4}

We are now ready to study the application of $U^{(i)}$ to the spaces $\mathcal{W}^{(i)}_{n,d}$. 
\begin{lemma}\label{lemmacrit}
For $m_1,m_2\in\mathbb{Z}$ and $n\in \mathbb{N}_0$,
\begin{align*}
U^{(1)}\left( \frac{\left(m_1x+m_2\right)r_L}{(2+11x)^n} \right)\in\mathcal{W}^{(0)}_{11n,n}.
\end{align*}
\end{lemma}
\begin{proof}
Analogous to the proof of Lemma \ref{auxiliarycongruence}. By \eqref{auxiliarylinefirst}, we need to verify the claim for $0\le n\le 32$.  Moreover, since $U^{(1)}$ linear, we only need to consider the cases $(m_1,m_2)\in\{(1,0),(0,1)\}$.  We verify these cases in Section 3a of the supplement file.
\end{proof}

Next, we study action of $U^{(1)}$ on $f$.

\begin{theorem}\label{u1tou0}
For $f\in\mathcal{W}^{(1)}_{n,0}$, we have
\begin{align*}
U^{(1)}\left( f \right)\in\mathcal{W}^{(0)}_{11n,n}.
\end{align*}
\end{theorem}

\begin{proof}
Let
\begin{align*}
f := \frac{1}{(2+11x)^{n}}\sum_{j+k\le 1}y^j z^k\sum_{m\ge 1-(j+k)} s_{jk}(m) 11^{\theta^{1}_{jk}(m)} x^m,
\end{align*} where $s_{jk}(m)\in\mathbb{Z}\left[\frac 12\right]$ has finite support in $m$. Note that
\begin{align*}
U^{(1)}(f) = \sum_{j+k\le 1}\sum_{m\ge 1-(j+k)} s_{jk}(m) 11^{\theta^{1}_{jk}(m)} U^{(1)}\left(\frac{y^j z^kx^m}{(2+11x)^{n}}\right).
\end{align*}  Now, by Theorem \ref{theoremcrit},
\begin{align*}
U^{(1)}(f) = \frac{1}{(2+11x)^{11n} c(x)^{n}} \sum_{\substack{j+k\le 1,\\ m\ge 1-(j+k)}}\sum_{\substack{0\le \gamma+\delta\le 1,\\ r\ge 0}} s_{jk}(m)\, h^{1jk}_{\gamma\delta}(m,n,r)11^{\pi^{1jk}_{\gamma\delta}(m,r)+\theta^{1}_{jk}(m)}y^{\gamma} z^{\delta} x^r.
\end{align*}  
Thus, it suffices to show that
\begin{align}
\pi^{1jk}_{\gamma\delta}(m,r)+\theta^{1}_{jk}(m)\ge \theta^{0}_{\gamma\delta}(r).\label{1to0ineq}
\end{align} 
There exists $M:=M(j,k)\in \mathbb{N}$ such that for $m\ge M$, we bound against
\begin{align*}
\pi^{1jk}_{\gamma\delta}(m,r)+\theta^{1}_{jk}(m) &\ge \left\lfloor \frac{11r-m+\lambda^{1jk}_{\gamma\delta}}{13} \right\rfloor + \left\lfloor \frac{11m+11j+22k-48}{13} \right\rfloor\ge \left\lfloor \frac{11r+10m+11j+22k+\lambda^{1jk}_{\gamma\delta}-60}{13} \right\rfloor\\
&\ge \left\lfloor \frac{11r+11\gamma+22\delta-7}{13} \right\rfloor + \left\lfloor \frac{10m+11j+22k-11\gamma-22\delta+\lambda^{1jk}_{\gamma\delta}-53}{13} \right\rfloor\\
&\ge \left\lfloor \frac{11r+11\gamma+22\delta-7}{13} \right\rfloor\ge \theta^{0}_{\gamma\delta}(r).
\end{align*} 
In particular, 
\[
M=\begin{cases}
6 &\text{ if }\ j=k=0,\\
5 &\text{ if }\ j=1,\\
4 &\text{ if }\ k=1.
\end{cases}
\]
For the remaining cases we check the inequalities more directly.  Below we check the inequalities for the floor values of
\begin{align*}
\left\lfloor\pi^{1jk}_{\gamma\delta}(m,r)+\theta^{1}_{jk}(m) - \theta^{0}_{\gamma\delta}(r)\right\rfloor.
\end{align*} 
Precisely, for $(j,k,\gamma,\delta)=(0,0,\gamma,\delta)$ and $1\le m\le 5$, we have
\begin{align*}
(0,0,0,0) &: \left\lfloor \frac{11r-m}{13} \right\rfloor - \left\lfloor \frac{11r-7}{13} \right\rfloor \ge 0,\\
(0,0,1,0) &: \left\lfloor \frac{11r+11-m}{13} \right\rfloor - \left\lfloor \frac{11r+4}{13} \right\rfloor \ge 0,\\
(0,0,0,1) &: \left\lfloor \frac{11r+24-m}{13} \right\rfloor - \left\lfloor \frac{11r+15}{13} \right\rfloor \ge 0.
\end{align*} 
Moreover, for $(j,k,\gamma,\delta)=(1,0,\gamma,\delta)$ and $0\le m\le 4$, it follows that
\begin{align*}
(1,0,0,0) &: \left\lfloor \frac{11r-m-1}{13} \right\rfloor - \left\lfloor \frac{11r-7}{13} \right\rfloor \ge 0,\\
(1,0,1,0) &: \left\lfloor \frac{11r+10-m}{13} \right\rfloor - \left\lfloor \frac{11r+4}{13} \right\rfloor \ge 0,\\
(1,0,0,1) &: \left\lfloor \frac{11r+23-m}{13} \right\rfloor - \left\lfloor \frac{11r+15}{13} \right\rfloor \ge 0.
\end{align*} 
Finally, for $(j,k,\gamma,\delta)=(0,1,\gamma,\delta)$ and $0\le m\le 3$, we obtain
\begin{align*}
(0,1,0,0) &: \left\lfloor \frac{11r-m-3}{13} \right\rfloor - \left\lfloor \frac{11r-7}{13} \right\rfloor \ge 0,\\
(0,1,1,0) &: \left\lfloor \frac{11r+8-m}{13} \right\rfloor - \left\lfloor \frac{11r+4}{13} \right\rfloor \ge 0,\\
(0,1,0,1) &: \left\lfloor \frac{11r+19-m}{13} \right\rfloor - \left\lfloor \frac{11r+15}{13} \right\rfloor \ge 0.
\end{align*}
This concludes the proof.
\end{proof}

\begin{theorem}\label{u0tou1}
For $g\in\mathcal{W}^{(0)}_{n,0}$ with $55\mid n$, there exist $m_1,m_2\in\mathbb{Z}$ such that
\begin{align*}
\frac{1}{11}U^{(0)}\left( g \right) - \frac{\left(m_1x+m_2\right)r_L}{(2+11x)^{11n+15}c(x)^n}\in\mathcal{W}^{(1)}_{11n+15,n}.
\end{align*}
\end{theorem}

\begin{proof}
Let
\begin{align*}
g := \frac{1}{(2+11x)^{n}}\sum_{j+k\le 1}y^j z^k\sum_{m\ge 1-(j+k)} s_{jk}(m) 11^{\theta^{0}_{jk}(m)} x^m,
\end{align*} with $s_{jk}(m)\in\mathbb{Z}\left[\frac 12\right]$ is with finite support.  We consider
\begin{align*}
U^{(0)}(g) = \sum_{j+k\le 1}\sum_{m\ge 1-(j+k)} s_{jk}(m) 11^{\theta^{0}_{jk}(m)} U^{(0)}\left(\frac{y^j z^kx^m}{(2+11x)^{n}}\right).
\end{align*}  Again applying Theorem \ref{theoremcrit},
\begin{align*}
U^{(0)}(g) = \frac{1}{(2+11x)^{11n+15} c(x)^{n}} \sum_{\substack{j+k\le 1,\\ m\ge 1-(j+k)}}\sum_{\substack{0\le \gamma+\delta\le 1,\\ r\ge 0}} s_{jk}(m)\, h^{0jk}_{\gamma\delta}(m,n,r)11^{\pi^{0jk}_{\gamma\delta}(m,r)+\theta^{0}_{jk}(m)}y^{\gamma} z^{\delta} x^r.
\end{align*}  

Next, we show that
\begin{align*}
\pi^{0jk}_{\gamma\delta}(m,r)+\theta^{0}_{jk}(m) \ge \theta^{1}_{\gamma\delta}(r)+1.
\end{align*}
First, for
\begin{equation}\label{ineq}
10m+11j+22k-11\gamma-22\delta+\lambda^{0jk}_{\gamma\delta}+16\ge 0,
\end{equation}
we have
\begin{align*}
\pi^{0jk}_{\gamma\delta}(m,r)+\theta^{0}_{jk}(m) &\ge \left\lfloor \frac{11r-m+\lambda^{0jk}_{\gamma\delta}}{13} \right\rfloor + \left\lfloor \frac{11m+11j+22k-7}{13} \right\rfloor\ge \left\lfloor \frac{11r+10m+11j+22k+\lambda^{0jk}_{\gamma\delta}-19}{13} \right\rfloor\\
&\ge \left\lfloor \frac{11r+11\gamma+22\delta-48}{13} \right\rfloor + \left\lfloor \frac{10m+11j+22k-11\gamma-22\delta+\lambda^{0jk}_{\gamma\delta}+29}{13} \right\rfloor\\
&\ge \left\lfloor \frac{11r+11\gamma+22\delta-48}{13} \right\rfloor + 1\ge \theta^{1}_{\gamma\delta}(r)+1,
\end{align*} 
Note that \eqref{ineq} holds for $m\ge 2$ if $j=k=0$, $m\ge 1$ if $j=1$ and $k=1$. Thus, for the remaining cases, we estimate lower bounds of
\begin{align}
\pi^{0jk}_{\gamma\delta}(m,r)+\theta^{0}_{jk}(m) - \theta^{1}_{\gamma\delta}(r) - 1.\label{importantinequalityA}
\end{align}

First, consider $j=k=0$. Since $m\ge 1-(j+k)$ and $r\ge 1-(\gamma+\delta)$, we need to verify for $m=1$. 
For $(j,k,\gamma,\delta)=(0,0,0,0)$ and $r\ge 5$,
\begin{align*}
\left\lfloor \frac{11r-35}{13} \right\rfloor - \left\lfloor \frac{11r-48}{13} \right\rfloor - 1\ge 0.
\end{align*}  
Next, for $(j,k,\gamma,\delta)=(0,0,1,0)$, and $r\ge 4$,
\begin{align*}
\left\lfloor \frac{11r-24}{13} \right\rfloor - \left\lfloor \frac{11r-37}{13} \right\rfloor - 1\ge 0.
\end{align*} 
Moreover, for $(j,k,\gamma,\delta)=(0,0,0,1)$ and $r\ge 3$,
\begin{align*}
\left\lfloor \frac{11r-13}{13} \right\rfloor - \left\lfloor \frac{11r-26}{13} \right\rfloor - 1\ge 0.
\end{align*} 

Now we consider $j=1$ and  we only need to verify for $m=0$. For $(j,k,\gamma,\delta)=(1,0,0,0)$ and $r\ge 5$,
\begin{align*}
\left\lfloor \frac{11r-34}{13} \right\rfloor - \left\lfloor \frac{11r-48}{13} \right\rfloor - 1.
\end{align*} 
For $(j,k,\gamma,\delta)=(1,0,1,0)$ and $r\ge 4$,
\begin{align*}
\left\lfloor \frac{11r-23}{13} \right\rfloor - \left\lfloor \frac{11r-37}{13} \right\rfloor - 1.
\end{align*} 
For $(j,k,\gamma,\delta)=(1,0,0,1)$ and $r\ge 3$,
\begin{align*}
\left\lfloor \frac{11r-12}{13} \right\rfloor - \left\lfloor \frac{11r-26}{13} \right\rfloor - 1.
\end{align*}  

Next, consider the case $k=1$ and \eqref{ineq} holds trivially. Moreover, we rewrite 
\begin{align*}
&\sum_{m\ge 0} s_{01}(m)11^{\theta^{(0)}_{01}(m)}6^m 4^n r_L=\sum_{m\ge 0} s_{01}(m)11^{\left\lfloor (11m+15)/13 \right\rfloor}6^m 4^n r_L = \sum_{m\ge 0} s_{01}(m)11^{\left\lfloor (11m+2)/13 \right\rfloor}6^m 4^n 11 r_L.
\end{align*}  
Since $s_{01}(m)$ has finite support, this is an integer multiple of $11r_L$.

It remains to analyze the negative powers of $11$ arising for the cases $j=k=0$ and $j=1$. With 
\[
1-(\gamma+\delta)\le r\le 4-(\gamma+2\delta),
\]
setting the coefficients of $y^{\gamma}z^{\delta}x^r$ as $t_{\gamma\delta}(r)$, we get
\begin{align}\label{t001a}
t_{\gamma\delta}(r) \equiv h^{000}_{\gamma\delta}(1,n,r)s_{00}(1)&+h^{010}_{\gamma\delta}(0,n,r)s_{10}(0)\pmod{11}.
\end{align}
We observe that only $s_{00}(1)$ and $s_{10}(0)$ contribute. Thus, we only consider the auxiliary coefficients for $U^{(0)}(\frac{x^jy^{(1-j)}}{(2+11x)^n})$ with $j\in\{0,1\}$. Since $55\mid n$, it is not hard to see that 
\begin{align}
h^{000}_{00}(1,11n,1)&\equiv 3\cdot 4^n\pmod{11}\label{h00000c}
\end{align} for $1\le n\le 3$.  
By Lemma \ref{auxiliarycongruence}, we see that
\begin{align*}
&\left\{h^{000}_{00}(1,11n,1)\!\pmod{11}\right\}_{n=1}^{\infty} = \{3,1,4,5,9,3,1,4,...\}.
\end{align*}  
Thus, \eqref{h00000c} is true for all $n\ge 0$ because the coefficients are periodic modulo $11$ with period $5$.  Consequently, for $55\mid n$, we have
\begin{align*}
h^{000}_{00}(1,n,1)&\equiv 9\pmod{11}.
\end{align*} 
So, \eqref{t001a} reduces to
\begin{align*}
t_{00}(1) &\equiv 9s_{00}(1)+5s_{10}(0)\pmod{11},\qquad 
t_{00}(2) \equiv 3s_{00}(1)+9s_{10}(0)\pmod{11},\\
t_{00}(3) &\equiv 9s_{00}(1)+5s_{10}(0)\pmod{11},\qquad 
t_{00}(4) \equiv 4s_{00}(1)+s_{10}(0)\pmod{11},\\
t_{10}(0) &\equiv 2s_{00}(1)+6s_{10}(0)\pmod{11},\qquad 
t_{10}(1) \equiv 4s_{00}(1)+s_{10}(0)\pmod{11},\\
t_{10}(2) &\equiv 9s_{00}(1)+5s_{10}(0)\pmod{11},\qquad
t_{10}(3) \equiv 7s_{00}(1)+10s_{10}(0)\pmod{11},\\
t_{01}(0) &\equiv 8s_{00}(1)+2s_{10}(0)\pmod{11},\qquad 
t_{01}(1) \equiv 3s_{00}(1)+9s_{10}(0)\pmod{11},\\
t_{01}(2) &\equiv 6s_{00}(1)+7s_{10}(0)\pmod{11}.
\end{align*}  
These yield that 
\begin{align*}
&\sum_{\substack{\gamma+\delta\le 1,\\1-(\gamma+\delta)\le r\le 4-(\gamma+2\delta)}}t_{\gamma\delta}(r)y^{\gamma}z^{\delta}x^r=\left(9 x + 3 x^2 + 9 x^3 + 4 x^4 + 2 y + 4 x y + 9 x^2 y + 7 x^3 y + 8 z + 
 3 x z + 6 x^2 z\right)s_{00}(1)\\ 
 &\hspace{5 cm}+ \left(5 x + 9 x^2 + 5 x^3 + x^4 + 6 y + x y + 5 x^2 y + 10 x^3 y + 2 z + 
 9 x z + 7 x^2 z\right)s_{10}(0)\\
&\hspace{1 cm}\equiv \left(1+x\right)\left(9s_{00}(1)+5s_{10}(0) \right)11 r_L\pmod{11},
\end{align*} and this concludes the proof.
\end{proof}

\subsection{Cancellation of $c(x)$}\label{cancelcx}

Let
\begin{align*}
\mathcal{F}^{(jk)} := \sum_{m\ge 0} s(m)11^{\theta^i_{jk}(m)} x^m.
\end{align*}  
Suppose that $c(x)^n|\mathcal{F}^{(jk)}(x)$, i.e., there exists a polynomial $\mathcal{F}^{(jk)}_0$ such that
\begin{align*}
\mathcal{F}^{(jk)}(x) = c(x)^n \mathcal{F}^{(jk)}_0.
\end{align*}  
In this section, our aim is to show that $\mathcal{F}^{(jk)}_0$ has a similar 11-adic form to that of $\mathcal{F}^{(jk)}$.
\begin{lemma}\label{factoringoutcx}
Let $M,m_1,m_2\in\mathbb{Z}$, with $\gcd(m_1m_2,11)=1$, and let $s(m)$ be an integer-valued function with finite support.  If
\begin{align*}
c(x)^n \Bigg| \left(i\left(m_1x+m_2\right)M r_L + \sum_{j+k\le 1} y^jz^k \sum_{m\ge 0} s(m)11^{\theta^i_{jk}(m)} x^m\right),
\end{align*} then there exists an integer $\hat{M}\in\mathbb{Z}$ and another integer-valued function $\hat{s}(m)$, also with finite support, such that
\begin{align*}
\frac{iM\left(m_1x+m_2\right) r_L + \sum_{j+k\le 1} y^jz^k \sum_{m\ge 0} s(m)11^{\theta^i_{jk}(m)} x^m}{c(x)^n} = i\hat{M}\left(m_1x+m_2\right) r_L + \sum_{j+k\le 1} y^jz^k \sum_{m\ge 0} \hat{s}(m)11^{\theta^i_{jk}(m)} x^m.
\end{align*}
\end{lemma}
\begin{proof}
Suppose that
\begin{align*}
\mathcal{F}_0^{(ij)}:= \sum_{m\ge 0} s_0(m)x^m.
\end{align*}  
We show that $s_0(m)$ is divisible by $11^{\theta^i_{jk}(m)}$ by induction principle. Since \[\left\lfloor\frac{11m_1+27}{13}\right\rfloor+\left\lfloor\frac{11m_2+27}{13}\right\rfloor\ge \left\lfloor\frac{11\left(m_1+m_2\right)+27}{13}\right\rfloor,\]
 $c(x)^n$ has the form
\begin{align*}
c(x)^n = C_n(0) + \sum_{r\ge 1} C_n(r)11^{\left\lfloor \frac{11r+27}{13} \right\rfloor}x^r.
\end{align*} 
Thus, we obtain
\begin{align*}
s(0)11^{\theta^i_{jk}(0)} &= C_n(0) s_0(0).
\end{align*} 
Since $C_n(0)=C_1(0)^n$ and $11\nmid C_n(0)$, there must exists $\hat{s}(0)\in \mathbb{Z}$ such that
\begin{align*}
s_0(0) = \hat{s}(0)11^{\theta^i_{jk}(0)}.
\end{align*}  This concludes the base case.  By induction hypothesis, $11^{\theta^i_{jk}(r)}\mid s_0(r)$ for $0\le r\le m-1$ and for some fixed $m\ge 1$.  Then for each $r$, there exists $\hat{s}(r)\in\mathbb{Z}$ such that
\begin{align*}
s_0(r) = \hat{s}(r)11^{\theta^i_{jk}(r)}.
\end{align*} 
Assume that, for $r\ge 4$,
\begin{align*}
\theta^i_{jk}(r) = \left\lfloor \frac{11r+\mu}{13} \right\rfloor,
\end{align*} with $\mu:=\mu(i,j,k)\in\mathbb{Z}$.  We then have
\begin{align}\label{ajkmfactorc4}
s(m)11^{\theta^i_{jk}(m)} &= \sum_{r=1}^{m-1} C_n(r) s_0(m-r)11^{\left\lfloor\frac{11r+27}{13}\right\rfloor} + C_n(0) s_0(m)\\
&\ge \sum_{r=1}^{m-1} C_n(r) \hat{s}(m-r)11^{\left\lfloor\frac{11r+14}{13}\right\rfloor+\left\lfloor\frac{11(m-r)+\mu}{13}\right\rfloor+1} + C_n(0) s_0(m)\nonumber\\
&\ge \sum_{r=1}^{m-1} C_n(r) \hat{s}(m-r) 11^{\left\lfloor\frac{11m+\mu+2}{13}\right\rfloor+1} + C_n(0) s_0(m)= 11^{\theta^i_{jk}(m)+1+\epsilon}\sum_{r=1}^{m-1} C_n(r) \hat{s}(m-r) + C_n(0) s_0(m)\nonumber
\end{align} with $\epsilon\ge 0$. Thus $11^{\theta^i_{jk}(m)}\mid s_0(m)$ and we write 
\begin{align*}
s_0(m) = \hat{s}(m)11^{\theta^i_{jk}(m)}.
\end{align*}  
Consequently, with $\hat{s}(m)\in\mathbb{Z}$, it follows that
\begin{align*}
\mathcal{F}^{(jk)}_0 = \sum_{m\ge 0} \hat{s}(m)11^{\theta^i_{jk}(m)}x^m.
\end{align*}

Now, we consider $i=1$ and $j=k=0$. With $t(m)\in\mathbb{Z}$ and $\gcd(t(m),11)=1$, we assume
\begin{align*}
s(0)&=0\ \ \text{and}\ \  s(m) = \frac{t(m)}{11}\ \text{for}\ 1\le m\le 3.
\end{align*}   
Following similar estimations done in \eqref{ajkmfactorc4}, for $1\le m\le 3$, we have
\begin{align*}
s_0(m)= \frac{\hat{s}(m)}{11}11^{\theta^1_{00}(m)}\ \ \text{and}\ \ \hat{s}(m)&\equiv C_n(0)^{-1}t(m)\pmod{11}.
\end{align*} 
Moreover, $s_0(m)$ is divisible by $11^{\theta^i_{00}(m)}$ for $m\ge 4$ and can be again show that following similar arguments in \eqref{ajkmfactorc4}. Analogous argument works as well for $i=j=1$, or for $i=k=1$. Let $m_1,m_2\in \mathbb{Z}$ with $\gcd(m_1m_2,11)=1$ and split $(m_1x+m_2)r_L$ as
\begin{align*}
\left(m_1x+m_2\right)r_L = r_L^{(00)}+yr_L^{(10)}+zr_L^{(01)}\ \ \text{with}\ \ r_L^{(jk)}\in\mathbb{Q}[x].
\end{align*}  Since $y$ and $z$ are not divisible by positive powers of $c(x)$, thus, if 
\begin{align*}
c(x)^n\mid  \left(m_1x+m_2\right)r_L + \sum_{j+k\le 1} y^jz^k \sum_{m\ge 0} s(m)11^{\theta^i_{jk}(m)} x^m=\sum_{j+k\le 1} y^jz^k\left(M\cdot r_L^{(jk)} + \sum_{m\ge 0} s(m)11^{\theta^i_{jk}(m)} x^m\right)
\end{align*}
we have
\begin{align*}
\frac{C_n(0)^{-1} M\left(m_1x+m_2\right)\left(x + 3 x^2 + 9 x^3 + 10 y + 10 x y + 2 x^2 y + 7 z + 8 x z\right)}{11}\equiv 0\pmod{\left(m_1+m_2x\right)r_L}.
\end{align*}
\end{proof}
\section{Proof of Theorem \ref{theorem0}}\label{completeproofA}

\begin{proof}[Proof of Theorem \ref{theorem0}]
By \eqref{L1identity11}, we have
\begin{align*}
\frac{2}{11} L_1 - \frac{(3+7x)r_L}{(2+11x)^{15}}\in\mathcal{W}^{(1)}_{15,0}.
\end{align*}  
From (\ref{psidefn}), $\psi(1)=15$. Now suppose that, for some $j_{\alpha},k_{\alpha}\in\mathbb{Z}\left[\frac 12\right]$, we have
\begin{align*}
\frac{1}{11^{\alpha}} L_{2\alpha-1} - \frac{(j_{2\alpha-1}+k_{2\alpha-1}x)r_L}{(2+11x)^{\psi(2\alpha-1)}} = f_{2\alpha-1}\in\mathcal{W}^{(1)}_{\psi(2\alpha-1),0}.
\end{align*}
Combining Lemma \ref{use}, Lemma \ref{lemmacrit}, and Theorem \ref{u1tou0}, we obtain
\begin{align*}
\frac{1}{11^{\alpha}} L_{2\alpha} &= \frac{1}{11^{\alpha}} U^{(1)}\left(L_{2\alpha-1}\right) = U^{(1)}\left( \frac{(j_{2\alpha-1}+k_{2\alpha-1}x)r_L}{(2+11x)^{\psi(2\alpha-1)}} \right) + U^{(1)}\left( f_{2\alpha-1} \right)\in\mathcal{W}^{(0)}_{11\psi(2\alpha-1),\psi(2\alpha-1)}.
\end{align*} 
It readily follows that 
\begin{align*}
11\psi(2\alpha-1) &= 11\left\lfloor \frac{11^{2\alpha}}{8} \right\rfloor = 11\left(\frac{11^{2\alpha}-1}{8}\right) = \frac{11^{2\alpha+1}-3}{8} - 1 = \psi(2\alpha).
\end{align*}
Thus, we get
\begin{align*}
\frac{1}{11^{\alpha}} L_{2\alpha} &=\mathcal{W}^{(0)}_{\psi(2\alpha),\psi(2\alpha-1)}.
\end{align*}  
Consequently, there exists
\begin{align*}
\mathcal{F}_{2\alpha}^{(jk)} = \sum_{m\ge 0} s_{jk}(m)11^{\theta^0_{jk}(m)} x^m\in\mathbb{Z}\left[\frac 12,x\right]\ \ \text{with}\ \ \frac{1}{11^{\alpha}} L_{2\alpha} =\frac{\mathcal{F}_{2\alpha}^{(00)}+y\mathcal{F}_{2\alpha}^{(10)}+z\mathcal{F}_{2\alpha}^{(01)}}{(2+11x)^{\psi(2\alpha)}c(x)^{\psi(2\alpha-1)}}.
\end{align*} 
So, by \eqref{Lalphainto0}, we have
\begin{align*}
\frac{(2+11x)^{\psi(2\alpha)}}{11^{\alpha}} L_{2\alpha} &=\frac{\mathcal{F}_{2\alpha}^{(00)}+y\mathcal{F}_{2\alpha}^{(10)}+z\mathcal{F}_{2\alpha}^{(01)}}{c(x)^{\psi(2\alpha-1)}}\in\bigoplus_{v\in\{x,y,z\}} v\,\mathbb{Z}\left[\frac 12,x\right],
\end{align*} whence
\begin{align*}
c(x)^{\psi(2\alpha-1)} | \mathcal{F}_{2\alpha}^{(jk)}.
\end{align*} 
Moreover, by Lemma \ref{factoringoutcx}, there exists
\begin{align*}
\hat{\mathcal{F}}_{2\alpha}^{(jk)} = \sum_{m\ge 0} \hat{s}_{jk}(m)11^{\theta^0_{jk}(m)} x^m\in\mathbb{Z}\left[\frac 12,x\right]\ \ \text{with}\ \ \frac{\mathcal{F}_{2\alpha}^{(jk)}}{c(x)^{\psi(2\alpha-1)}} = \hat{\mathcal{F}}_{2\alpha}^{(jk)}.
\end{align*}
Thus, we get
\begin{align*}
\frac{(2+11x)^{\psi(2\alpha)}}{11^{\alpha}} L_{2\alpha} &=\hat{\mathcal{F}}_{2\alpha}^{(00)}+y\hat{\mathcal{F}}_{2\alpha}^{(10)}+z\hat{\mathcal{F}}_{2\alpha}^{(01)}\in\bigoplus_{v\in\{x,y,z\}} v\,\mathbb{Z}\left[\frac 12,x\right]\ \ \text{and}\ \ \frac{1}{11^{\alpha}} L_{2\alpha}=\mathcal{W}^{(0)}_{\psi(2\alpha),0}.
\end{align*}

Next, applying $U^{(0)}$ and by Lemma \ref{u0tou1}, for some $\hat{j}_{\alpha+1},\hat{k}_{\alpha+1}\in\mathbb{Z}\left[\frac 12\right]$, we have
\begin{align*}
\frac{1}{11^{\alpha}} L_{2\alpha+1} &= \frac{1}{11^{\alpha}} U^{(0)}\left(L_{2\alpha}\right)=11 f_{2\alpha+1} + \frac{11 (\hat{j}_{2\alpha+1}+\hat{k}_{2\alpha+1}x)r_L}{(2+11x)^{11\psi(2\alpha)+15}c(x)^{\psi(2\alpha)}}\ \ \text{with}\ \ f_{2\alpha+1}\in\mathcal{W}^{(0)}_{11\psi(2\alpha)+15,\psi(2\alpha)}.
\end{align*}
Note that
\begin{align*}
11\psi(2\alpha)+15 = \psi(2\alpha+1),
\end{align*} and so,
\begin{align*}
\frac{1}{11^{\alpha+1}} L_{2\alpha+1} - \frac{(\hat{j}_{2\alpha+1}+\hat{k}_{2\alpha+1}x)r_L}{(2+11x)^{\psi(2\alpha+1)}c(x)^{\psi(2\alpha)}} = f_{2\alpha+1} \in\mathcal{W}^{(1)}_{\psi(2\alpha+1),\psi(2\alpha)}.
\end{align*} 
By \eqref{Lalphainto0}, $(2+11x)^{\psi(2\alpha+1)}$ cancels the pole of $L_{2\alpha+1}$ at $[\frac 12]$, the numerator of 
\begin{align*}
f_{2\alpha+1} + \frac{(\hat{j}_{2\alpha+1}+\hat{k}_{2\alpha+1}x)r_L}{(2+11x)^{\psi(2\alpha+1)}c(x)^{\psi(2\alpha)}}
\end{align*} is divisible by the factor $c(x)^{\psi(2\alpha)}$. So, by Lemma \ref{factoringoutcx}, there exist $j_{\alpha+1},k_{\alpha+1}\in\mathbb{Z}\left[\frac 12\right]$ such that
\begin{align*}
\frac{1}{11^{\alpha+1}} L_{2\alpha+1} - \frac{(j_{2\alpha+1}+k_{2\alpha+1}x)r_L}{(2+11x)^{\psi(2\alpha+1)}}\in\mathcal{W}^{(1)}_{\psi(2\alpha+1),0}.
\end{align*}
This concludes the proof.
\end{proof}

\section{Proof of Initial Relations for Theorem \ref{theoremcrit}}\label{proofoffundamentalrelations}

The proofs in Sections 2--4 depend on the following finitely many identities between modular functions, all of which we verify by computer:
\begin{enumerate}
	\item the modular equation \eqref{modX} for $x$, and hence \eqref{modT} for $t=2+11x$,
	\item the representation \eqref{L1identity11} of $L_1$,
	\item the $198$ relations \eqref{fundrelate} for $U^{(i)}(y^jz^kx^m)$, $i,j,k\in\{0,1\}$, $j+k\le1$, $0\le m\le32$, from which we obtain the $6534$ initial cases $0\le m,n\le32$ of Theorem \ref{theoremcrit},
	\item the initial cases $0\le n\le32$ of Lemma \ref{lemmacrit}.
\end{enumerate}
We now explain how each of these is reduced to a finite computation. Precisely, to complete the proofs, we need to confirm a certain set of initial relations between modular functions. Note that by \eqref{modT}, we have 
\begin{align*}
&U^{(i)}\left(\frac{y^jz^kx^m}{(2+11x)^n}\right) = -\frac{1}{b_{0}}\sum_{\nu=1}^{33} b_{\nu} U^{(i)}\left(\frac{y^jz^kx^m}{(2+11x)^{n-{\nu}}}\right) =\frac{2^{20}}{(2+11x)^{11}c(x)}\sum_{\nu=1}^{33} b_{\nu} U^{(i)}\left(\frac{y^jz^kx^m}{(2+11x)^{n-{\nu}}}\right).
\end{align*}  We can of course extrapolate until $n-\nu<0$, i.e., until we are studying combinations of
\begin{align*}
U^{(i)}\left(y^jz^kx^m(2+11x)^n\right)
\end{align*} for $n> 0$.  This in turn can be expressed through the linearity of $U^{(i)}$ and the binomial expansion $(2+11x)^n$ into combinations of 
\begin{align}
U^{(i)}\left(y^jz^kx^m\right).\label{fundrelate}
\end{align}  Thus, if we have explicit relations for \eqref{fundrelate} for $i,j,k\in\{0,1\}$ and $0\le m\le 32$, then from these we can construct the $6534$ relations for $U^{(i)}\left(\frac{y^jz^kx^m}{(2+11x)^n}\right)$ with $0\le m,n\le 32$.

Likewise, the proof of \eqref{modX} constitutes a polynomial relationship between modular functions that needs to be checked.

Relations for \eqref{fundrelate} recorded in the Mathematica supplement have the following form:
\begin{align}
U^{(1)}\left(y^jz^kx^m\right) &= f^{(1)}_{j,k,m}\in\bigoplus_{v\in\{x,y,z\}} v\,\mathbb{C}[x],\label{variousexamples0}\\
(2+11x)^{15}U^{(0)}\left(y^jz^kx^m\right) &= f^{(0)}_{j,k,m}\in\bigoplus_{v\in\{x,y,z\}} v\,\mathbb{C}[x],\label{variousexamples}
\end{align} in which $f^{(i)}_{j,k,m}$ are polynomials, given explicitly.  In principle, we could compare the principal parts of $(2+11x)^{15(1-i)}U^{(i)}\left(y^jz^kx^m\right)$ and $f^{(i)}_{j,k,m}$ at $[0]$.  However, this is difficult to compute, especially on the left-hand side of the relations.  On the other hand, we can multiply both sides of each relation by an eta quotient which annihilates the pole at $[0]$ and induces a pole at $[\infty]$.  The resulting principal parts are much easier to check. Indeed, $\frac 1z$ suffices for our purpose.  There exists $M\in \mathbb{N}$ such that
\begin{align*}
U^{(i)}\left(\frac{\left(2+11x(11\tau)\right)^{15(1-i)}}{z(11\tau)^M}y^jz^kx^m\right),\ \frac{1}{z^M} f^{(i)}_{j,k,m}\in\mathcal{M}^{\infty}\left(  \mathrm{X}_0(22)\right).
\end{align*}  Remembering the definition of $U^{(i)}$, we need only check that the principal part of 
\begin{align*}
U_{11}\left(\mathcal{A}^{1-i}\frac{\left(2+11x(11\tau)\right)^{15(1-i)}}{z(11\tau)^M}y^jz^kx^m\right) - \frac{1}{z^M} f^{(i)}_{j,k,m}
\end{align*} vanishes. These can be computed in Section 5 of the Mathematica supplement.  Note that $L_1$ constitutes the case of (\ref{variousexamples}) for $j=k=m=0$.  The modular equation \eqref{modX} can likewise be computed by cusp analysis in the notebook.  The calculation of (\ref{variousexamples0})-(\ref{variousexamples}) and \eqref{modX} is done with PARI/GP.\footnote{This is because computing with Mathematica takes substantial time.} 

\section{Conclusion}\label{sectconcl}

We conclude this paper by proposing the following conjecture.
For $a\in\mathbb{Z}$ define $p^{[1]}_a(n)$ and $p^{[2]}_a(n)$ by
\[
P^{[1]}_a(q):=\frac{(-q;q)_\infty^a}{(q;q)_\infty}=:\sum_{n\ge0}p^{[1]}_a(n)q^n,
\qquad
P^{[2]}_a(q):=\frac{(-q;q)_\infty^a}{(q;q)_\infty^5}=:\sum_{n\ge0}p^{[2]}_a(n)q^n .
\]
Since $(-q;q)_\infty=(q^2;q^2)_\infty/(q;q)_\infty$, these are eta quotients of level $2$, and the associated congruence families live on $\mathrm{X}_0(22)$, the very curve of genus $2$ that we have hitherto studied. Note that $P^{[2]}_2(q)=D_2(q)$, so Theorem \ref{Thm12} is the case $a=2$ of part (2) below.

\begin{conjecture}\label{conj:families}
Let $\alpha \in \mathbb{N}$ and $n \in \mathbb{N}_0$.
\begin{enumerate}
	\item If $a \equiv 3, 5, 7 \pmod{11}$ and $|a| \le 35$, then
	\[
	p^{[1]}_a\left(11^{\alpha} n + \delta^{[1]}_{a,\alpha}\right) \equiv 0\ \left(\textnormal{mod}\ 11^{\lceil \alpha/2 \rceil}\right),
	\quad\text{where } 24\,\delta^{[1]}_{a,\alpha} \equiv 1 - a \pmod{11^{\alpha}}.
	\]
	\item If $a \not\equiv 1, 3, 5 \pmod{11}$ and $-44 \le a \le 39$, then
	\[
	p^{[2]}_a\left(11^{\alpha} n + \delta^{[2]}_{a,\alpha}\right) \equiv 0\ \left(\textnormal{mod}\ 11^{\lceil \alpha/2 \rceil}\right),
	\quad\text{where } 24\,\delta^{[2]}_{a,\alpha} \equiv 5 - a \pmod{11^{\alpha}}.
	\]
\end{enumerate}
	\end{conjecture}
It is natural to ask which methods could prove Conjecture \ref{conj:families}. The classical approach via modular equations, in the tradition of Watson and Atkin and as refined by Paule--Radu and by the localization method, does succeed, as the present paper shows, but only at the cost of a non-unit localizing factor, an auxiliary denominator $c(x)$ that must be properly controlled, and a substantial computation. A more conceptual argument would clearly be highly desirable. The $\ell$-adic Galois-theoretic methods developed by Folsom--Kent--Ono, Boylan--Webb guarantee that congruences of this shape exist. But it is uncertain whether the exponent $\lceil\frac{\alpha}2\rceil$ can be properly accounted for. Moreover, the generating functions involved are weakly holomorphic forms of half-integral weight with an eta multiplier, not Hecke eigenforms to which such techniques most naturally apply. The framework developed in \cite{BGKMR} based on slopes of $U_\ell$ on over-convergent forms is perhaps the most suggestive, since the growth rate $\lceil\frac{\alpha}2\rceil$ is strongly reminiscent of a slope $\frac12$. However, this would require packaging these forms in the appropriate component of weight space at a certain tame level.

\section{Appendix}\label{Appen}

Here we state explicitly the expression for $c(x)$.
\begin{align*}
c(x)&:= 69550975755977284911104\!-\!71287715497527282651820032 x\\
&+ 17900193598896240733994540032 x^2 +2433656434135752614962897574400 x^3\\
&+ 135334237711138006616618300344576 x^4 + 4603551043973265247618043536808704 x^5\\
&+ 109616062492953370639775636433454464 x^6 + 1954419189636606079380190879818816064 x^7\\
&+ 27144272314898816695853907223593334304 x^8 + 301089044785364403767185228909074692504 x^9\\
&+ 2710730798554339991175805693674009005740 x^{10}+ 20011291470089564255540586248512208139082 x^{11}\\
&+ 121821501325881013307089187765394122218837 x^{12}+ 612703095907075970646910882428597137091912 x^{13}\\
&+ 2541715085550544936276941583396435763817408 x^{14}+ 8649797863778836226399665129987964376904896 x^{15}\\
&+ 23913154138234188165042955394025759141437184 x^{16}+ 52880092072078116046416414506583439275205632 x^{17}\\
&+ 91342853380107788221814612586966899809652736 x^{18}+ 118783323283193765986323689131241918352015360 x^{19}\\
&+ 109398375434451346894431470218757768960606208 x^{20}+ 63625268819176226667199093217678019620241408 x^{21}\\
&+ 17574243049798037033395145952709018207125504 x^{22}.
\end{align*}

Here are the initial values of the constants $\lambda^{ijk}_{\gamma\delta}$ from Theorem \ref{theoremcrit}.

\begin{align*}
&\lambda^{100}_{00} = 0,\ \ \ \lambda^{000}_{00} = -34\ \ \ \lambda^{100}_{10} = 11,\ \ \ \lambda^{000}_{10} = -23\ \ \ \lambda^{100}_{01} = 24,\ \ \ \lambda^{000}_{01} = -12,\\
&\lambda^{110}_{00} = -1,\ \ \ \lambda^{010}_{00} = -34\ \ \ \lambda^{110}_{10} = 10,\ \ \ \lambda^{010}_{10} = -23\ \ \ \lambda^{110}_{01} = 23,\ \ \ \lambda^{010}_{01} = -12,\\
&\lambda^{101}_{00} = -3,\ \ \ \lambda^{001}_{00} = -38\ \ \ \lambda^{101}_{10} = 8,\ \ \ \lambda^{001}_{10} = -27\ \ \ \lambda^{101}_{01} = 19,\ \ \ \lambda^{001}_{01} = -16.
\end{align*}

\section{Acknowledgments}
This research began when the first author was funded by the Austrian Science Fund (FWF): W1214-N15, Project DK6 and currently, he is supported by the European Research Council (ERC) under the European Unions Horizon 2020 research and innovation programme (grant agreement No. 101001179). Nicolas Allen Smoot was funded in whole by the Austrian Science Fund (FWF) Principal Investigator Project (No.~10.55776/ PAT2615425, ``Multiplicities Between Modular Congruences"). He would like to thank the Austrian Government and People for their generous support.  A personal thanks to Denise and Alexander Teuchtmann for their inspiration during the writing of this paper.


\begin{thebibliography}{X}
	
\bibitem{Ahlgren2000} S.~Ahlgren, \emph{Distribution of the partition function modulo composite integers $M$}, Math. Ann. \textbf{318} (2000), 795--803.

\bibitem{AhlgrenAllenTang2022} S.~Ahlgren, P.~B.~Allen, and S.~Tang, \emph{Congruences like Atkin's for the partition function}, Trans. Amer. Math. Soc. Ser. B \textbf{9} (2022), 1044--1064.

\bibitem{AhlgrenBoylan2003} S.~Ahlgren and M.~Boylan, \emph{Arithmetic properties of the partition function}, Invent. Math. \textbf{153} (2003), 487--502.

\bibitem{AhlgrenOno2001} S.~Ahlgren and K.~Ono, \emph{Congruence properties for the partition function}, Proc. Natl. Acad. Sci. USA \textbf{98} (2001), 12882--12884.




	
	


\bibitem{AndrewsPaule} G. Andrews and P. Paule, {\it MacMahon's partition analysis XIII: Schmidt type partitions and modular forms}, J. Number Theory. {\bf 234} (2022), 95--119.

\bibitem{Atkin} A. Atkin, {\it Proof of a conjecture of Ramanujan}, Glasgow Math. J. {\bf 8} (1967), 14--32.


\bibitem{Banerjee1} K. Banerjee and N. A. Smoot, {\it The localization method applied to $k$-elongated plane partitions and divisibility by $5$}, Math. Z. {\bf 309} (46), (2025).

\bibitem{Banerjee2} K. Banerjee and N. A. Smoot, {\it $2$-Elongated plane partitions and powers of $7$: the localization method applied to a genus $1$ congruence family}, Canad. J. Math. (2025), 1--43.



\bibitem{Beazer} M. Beazer, \textit{3-adic properties of Hecke traces of singular moduli}, Master thesis, Brigham Young University (2021), online available at \url{http://hdl.lib.byu.edu/1877/etd11831}.

\bibitem{BW} M.~Boylan and J.~Webb, {\it The partition function modulo prime powers}, Trans. Amer. Math. Soc. \textbf{365} (2013), 2169--2206.

\bibitem{BGKMR} K.~Bringmann, P.~Guerzhoy, B.~Kane, M.~Mertens, and L.~Rolen, {\it On $U_p$-congruences for meromorphic modular forms via supersingularity}, (2026), arXiv:2606.14020.




\bibitem{Calegari} F. Calegari, {\it Congruences between modular forms}, Arizona winter school (2013), online available at \url{https://www.math.uchicago.edu/~fcale/papers/AWS.pdf}.



\bibitem{FKO} A.~Folsom, Z.~Kent, and K.~Ono, \emph{$\ell$-adic properties of the partition function} (with an appendix by N.~Ramsey), Adv. Math. \textbf{229} (2012), no.~3, 1586--1609.



\bibitem{Forster} O. Forster, \textit{Lectures on Riemann Surfaces,} Springer-Verlag (1981).

\bibitem{Garvan0} F. Garvan, J. Sellers, and N. A. Smoot, {\it Old meets new: connecting two infinite families of congruences modulo powers of $5$ for generalized Frobenius partition functions}, Adv. Math. \textbf{454} (2024), p. 109866.

\bibitem{Gordon} B. Gordon and K. Hughes, {\it Ramanujan Congruences for $q(n)$},  In: M. Knopp (eds) analytic number theory, lecture notes in mathematics {\bf  899}, Springer, Berlin, Heidelberg (1981). 

\bibitem{Guerzhoy} P. Guerzhoy, {\it Congruences for traces of singular moduli and Hurwitz--Kronecker class numbers}, (2026), arXiv:2602.19544.

\bibitem{Ono2000} K.~Ono, \emph{Distribution of the partition function modulo $m$}, Ann. of Math. (2) \textbf{151} (2000), 293--307.	





\bibitem{Paule}  P. Paule and S. Radu, {\it The Andrews--Sellers family of partition congruences}, Adv. Math. {\bf 230} (2012), 819--838.

\bibitem{Paule2} P. Paule and S. Radu, {\it A Proof of the Weierstrass Gap Theorem Not Using the Riemann--Roch Formula}, Ann. Comb. {\bf 23}, 963-1007 (2019).

\bibitem{Radu} S. Radu, {\it An algorithmic approach to Ramanujan--Kolberg identities}, J. Symb. Comp. {\bf 68} (2015), 225--253.


\bibitem{Ramanujan} S. Ramanujan, {\it Some properties of $p(n)$, the number of partitions of $n$}, Proc. Cambridge Phil. Soc. {\bf 19} (1919), 207--210.



\bibitem{Smoot2} N. A. Smoot, ``A Congruence Family For 2-Elongated Plane Partitions: An Application of the Localization Method", \textit{Journal of Number Theory} 242, pp. 112--153 (2023).



\bibitem{Watson} G. Watson, {\it Ramanujans Vermutung über Zerfällungsanzahlen}, J. Reine Angew. Math. {\bf 179} (1938), 97--128.

\end{thebibliography}
\end{document}